\documentclass[11pt]{article}
\usepackage[utf8]{inputenc}
\usepackage{caption}
\usepackage{subcaption}
\usepackage[margin=1in]{geometry}
\usepackage{amsmath}
\usepackage{amsfonts}
\usepackage{graphicx}
\usepackage{float}
\usepackage{listings}
\usepackage{hyperref}
\usepackage{matmacros}

\author{Cody J. Balos\\
\small{Lawrence Livermore National Laboratory}}
\title{Data-driven time-scale separation of ODE right-hand sides using dynamic mode decomposition and time delay embedding}
\date{}

\begin{document}

\maketitle
\abstract{Multi-physics simulation often involve multiple different scales. The ARKODE
ODE solver package in the SUNDIALS library addresses multi-scale problems with a multi-rate time-integrator that
can work with a right-hand side that has fast scale and slow scale components. In this report, we use
dynamic mode decomposition and time delay embedding to extract the fast and and slow components of the
right-hand sides of a simple ODE from data. We then use the extracted components to solve the ODE with ARKODE.
Finally, to move towards a real-world use case, we attempt to extract fast and slow scale dynamics
from synthetic seismic modeling data.}

\begin{center}
    \noindent\rule{5cm}{0.4pt}~~$\clubsuit$~~\rule{5cm}{0.4pt}
\end{center}

\section{Introduction and Overview}

As computing power increases, multi-physics simulations which involve
the coupling of different physical phenomena are becoming more common and more ambitious.
The different physics in these simulations often occur at different scales. As a result,
there are inefficiencies when solving the ODE systems that arise because the time-step
of the integrator is restricted by the fast scale. The ARKODE \cite{reynoldsARKODE2023} ODE solver package in the SUNDIALS library \cite{gardnerSUNDIALS2022}
addresses these problems with a multirate time-integrator that is defined by a
right-hand side composed of a fast-scale and slow-scale components as shown in 
Equation~(\ref{eq:multirate}).
Naturally, this leads to the question of how to split the right-hand side appropriately.
In this report, we will take a data-driven approach based on the dynamic mode
decomposition (DMD), a powerful tool for approximating dynamics from data
\cite{doi:10.1137/1.9781611974508.ch1,2013arXiv1312.0041T}.

\begin{equation}\label{eq:multirate}
    \dot{y} = \underbrace{f_s(t,y)}_{\text{slow scale}} 
            + \underbrace{f_f(t,y)}_{\text{fast scale}},\quad y(t_0) = y_0
\end{equation}

In \cite{grosekDMDVideo2014}, Grosek \textit{et al.}, introduced
a method for extracting the low-rank and sparse components of a DMD approximation
in the context of background and foreground separation in videos, however, we will
apply the concept to other dynamical systems where there exists a slow-scale
(background) that is much slower than a fast-scale (foreground). In order to successfully
reconstruct the dynamics from systems with only one partial state information, we will use
time delay embeddings to augment the state. This method was presented by Kamb \textit{et al.}
as a Koopman observable in \cite{kamb2020time}.

After extracting the the fast and slow components of an simple ODE right-hand side, we
will solve a multi-rate ODE by feeding the approximated components into ARKODE. Finally,
we will attempt to extract fast and slow-scale dynamics from synthetic tsunami modeling data
produced with a code developed by Vogl and Leveque based on the Clawpack solver package
\cite{vogl2017high}. The idea is that with some additional work, the extracted components
could be used as the fast and slow right hand sides in ARKODE to solve ODE system(s) that stem from
the wave propagation equations discussed.

The rest of this report is organized as follows. Section~\ref{s:theory} discusses
the theory behind the methods used to form our separation algorithm.
Section~\ref{s:implement} explains the algorithm and its implementation in MATLAB.
Section~\ref{s:results} discussed the computational results for the toy problem
and the synthetic tsunami data. Finally, we provide our conclusions in Section~\ref{s:summary}. 

\section{Theoretical Background}\label{s:theory}

In this section we will discuss the theory behind the different methods used in our separation
algorithm. Consider a dynamical system

\begin{equation}\label{eq:dynamical}
    \frac{d\mathbf{x}}{dt} = \mathbf{f}(\mathbf{x},t;\mu),
\end{equation}

\noindent where $\mathbf{x}(t) \in \mathbb{R}^n$ is a vector representing the state of our dynamical system at time $t$, $\mu$ contains the parameters of the system, and $\mathbf{f}(\cdot)$ represents the dynamics. By sampling (\ref{eq:dynamical}) with a uniform sampling rate we get the discrete-time flow-map

\begin{equation}\label{eq:flow}
    \mathbf{x}_{k+1} = \mathbf{F}(\mathbf{x}_k).
\end{equation}

\subsection{Koopman Theory}

The Koopman operator $\mathcal{K}_t$ provides an operator perspective on dynamical systems.
The insight behind the operator is that the finite-dimensional, nonlinear dynamics
of (\ref{eq:dynamical}) can be transformed to an infinite-dimensional, linear dynamical
system by considering a new space of scalar observable functions 
$g: \mathbb{R}^N \rightarrow \mathbb{C}$ on the state instead of the state directly
\cite{kamb2020time}. The Kooopman operator acts on the observable:

\[
    \mathcal{K}_t g(\mathbf{x}_0) = g(\mathbf{F}_t(\mathbf{x}_0)) = g(x(t))  
\]

\noindent where $\mathcal{K}_t$ maps the measurement function $g$ to the next set of 
value it will take after time $t$.

\subsection{Dynamic Mode Decomposition}

By sampling (\ref{eq:dynamical}) with a uniform sampling rate, the DMD approximates the low-dimensional modes of the linear, time-independent Koopman operator in order to estimate the dynamics of the system. Consider $n$ data points with a total of $m$ samplings in time, then we 
can form the two matrices:

\begin{equation}\label{eq:X1}
    \mathbf{X_1} =
    \begin{bmatrix}
        \mid & \mid & \mid \\
        \mathbf{x_1} & \mathbf{x_2} & \cdots & \mathbf{x}_{m-1}\\
        \mid & \mid & \mid
    \end{bmatrix}
\end{equation}

\begin{equation}\label{eq:X2}
    \mathbf{X_2} =
    \begin{bmatrix}
        \mid & \mid & \mid \\
        \mathbf{x_2} & \mathbf{x_3} & \cdots & \mathbf{x}_m\\
        \mid & \mid & \mid
    \end{bmatrix}.
\end{equation}

The Koopman operator $\mathbf{A}$ provides a mapping that takes the system, modeled by
the data, from time $j$ to $j+1$, i.e. $x_{j+1} = \mathbf{A}x_j$. We can use the
singular value decomposition (SVD) of the matrix $\mathbf{X}_1$ to reduce the
rank $\ell$ as well as to obtain the benefits of unitary matrices. Accordingly,
$\mathbf{X_1} = \mathbf{U}\mathbf{\Sigma}\mathbf{V}^*$, where 
$\mathbf{U} \in \mathbb{C}^{n \times \ell}$ is unitary, $\mathbf{\Sigma} \in \mathbb{C}^{\ell\times\ell}$
is diagonal, and $\mathbf{U} \in \mathbb{C}^{(m-1) \times \ell}$ is unitary. Thus, the approximate 
eigendecomposition of the Koopman operator is:

\[
    \tilde{\mathbf{A}} = \mathbf{U}^*\mathbf{X_2}\mathbf{V}\mathbf{\Sigma},
\]

\noindent the eigenvalues are given by $\tilde{\mathbf{A}}\omega_j = \lambda_j\omega_j$, and
the DMD mode (eigenvector of the Koopman operator) corresponding to eigenvalue $\lambda_j$ is ${\varphi}_j = \mathbf{U}\omega_j$.

Using the approximate eigendecomposition of the Koopman operator we obtain the DMD
reconstruction of the dynamics:

\begin{equation}\label{eq:xdmd}
    \mathbf{x}_{\text{DMD}}(t) = \sum_{j=1}^\ell b_j \varphi_j e^{\omega_j t} = \mathbf{\Phi}\mathbf{\Omega^t}\mathbf{b},
\end{equation}

\noindent where

\[
    \mathbf{\Omega}=
    \begin{bmatrix}
        e^{\omega_1} & 0 & \cdots & 0 \\
        0 & e^{\omega_2} & \ddots & 0 \\
        \vdots & \ddots & \ddots & \vdots \\
        0 & 0 & \cdots & e^{\omega_t}
    \end{bmatrix}
\]

\noindent the vector $\mathbf{b} \approx \mathbf{\Phi}^{-1}\mathbf{x_1}$, and $\Phi = \mathbf{U}\mathbf{\Omega}$ \cite{grosekDMDVideo2014,2013arXiv1312.0041T}.

If we consider components of a dynamical system that changes very slowly with time with
respect to the fast components, it becomes clear that they will have an associated
 mode $\omega_j$ such that

\begin{equation}\label{eq:omega}
    \|\omega_j\|/\text{max}(\|\omega_j\|) \approx 0.
\end{equation}

For this reason the DMD can be used to separate the dynamics into a fast and slow component. Assume
that $\omega_{p}$, where $p \in \{1,2,\ldots,\ell\}$, satisfies Equation~(\ref{eq:omega}), and that 
$\|\omega_{j}\|/\text{max}(\|\omega_j\|) \ \forall \ j \neq p$ is bounded away from zero. Then we can redefine Equation~(\ref{eq:xdmd})
for a multi-scale dynamical system as:

\begin{equation}\label{eq:slowfast}
	{\bf X}_{\text{DMD}} = \underbrace{b_{p}\mathbf{\varphi}_{p}e^{\omega_{p} {\bf t}}}_{\text{slow scale}} + \underbrace{\sum_{j \neq p} b_{j}\mathbf{\varphi}_{j}e^{\omega_{j} {\bf t}}}_{\text{fast scale}}.
\end{equation}

If we consider the low-rank reconstruction (\ref{eq:lowrank}) of the DMD and
the fact that (\ref{eq:true}) should be true, it follows that the sparse reconstruction given can be calculated
using the real values of elements only.

\begin{equation}\label{eq:lowrank}
	{\bf X}_{\text{DMD}}^{\text{Low-Rank}} = b_{p}\mathbf{\varphi}_{p}e^{\omega_{p} {\bf t}},
\end{equation}

\begin{equation}\label{eq:true}
	{\bf X} = {\bf X}_{\text{DMD}}^{\text{Low-Rank}} + {\bf X}_{\text{DMD}}^{\text{Sparse}}.
\end{equation}

\begin{equation}
    {\bf X}_{\text{DMD}}^{\text{Sparse}} = \sum_{j \neq p} b_{j}\mathbf{\varphi}_{j}e^{\omega_{j} {\bf t}}    
\end{equation}

\begin{equation}\label{eq:sparse2}
	{\bf X}_{\text{DMD}}^{\text{Sparse}} = {\bf X} - \Big|{\bf X}_{\text{DMD}}^{\text{Low-Rank}}\Big|.
\end{equation}

\subsection{Time Delay Embedding}
\textit{Delay embedding} is a classic technique for overcoming partial state information.
In the context of the DMD and Koopman operator theory, we can use time delay embedding
to construct a new observable: $\tilde{\mathbf{g}}(\mathbf{x}(t)):=(g(\mathbf{x}(t)),g(\mathbf{x}(t-\Delta t)),g(\mathbf{x}(t-2\Delta t)),
g(\mathbf{x}(t-n\Delta t))) \in \mathbb{R}^n$ with lag time $\Delta t$. This is the delay embedding of the trajectory $g(\mathbf{x}(t))$. We can use the embedding to form the Hankel matrix:

\begin{equation}\label{eq:hankel}
    \mathbf{H}=
    \begin{bmatrix}
        g(\mathbf{x}_1) & g(\mathbf{x}_2) & \cdots & g(\mathbf{x}_M) \\
        g(\mathbf{x}_2) & g(\mathbf{x}_3) & \cdots & g(\mathbf{x}_{M+1}) \\
        \vdots & \vdots & \ddots & \vdots \\
        g(\mathbf{x}_N) & g(\mathbf{x}_N) & \cdots & g(\mathbf{x}_{N+M+1}) \\
    \end{bmatrix}
    =
    \begin{bmatrix}
        \mid & \mid & & \mid \\
        \tilde{\mathbf{g}}(\mathbf{x_1}) & \tilde{\mathbf{g}}(\mathbf{x_2}) & \cdots & \tilde{\mathbf{g}}(\mathbf{x_M}) \\
        \mid & \mid & & \mid
    \end{bmatrix}
\end{equation}

By examining the singular values of the Hankel matrix we can find the best rank-$r$ approximation of the trajectory space. Then, we can form the matrices $\mathbf{X_1}$ and $\mathbf{X_2}$ for the
DMD from $\mathbf{H}$ and set $\ell = r$ for the best DMD approximation using the embedding
\cite{kamb2020time}.

\section{Algorithm Implementation and Development}\label{s:implement}

We develop the algorithm with the assumption we are given $m$ measurements and $n$ snapshots
of a dynamical system with fast and slow scales, but only partial state information. This data forms the matrix $\mathbf{X} \in \mathbb{R}^{m \times n}$. Since we only have partial state
information, the first step in the algorithm is to employ time delay embeddings to form
a Hankel matrix $\mathbf{H}$ with $N$ embeddings. We choose an $N$ that minimizes the
approximation error in the second stage of the algorithm.

The second step in the algorithm is to form the matrices $\mathbf{X_1}$ and $\mathbf{X_2}$
from $\mathbf{H}$, and then apply the DMD method with a rank-reduction that is
chosen based on the singular values of $\mathbf{H}$. This produces the DMD approximation to the
dynamical system as given in Equation~(\ref{eq:xdmd}).

The third step in the algorithm is to extract the fast and slow components by applying
Equations~(\ref{eq:slowfast}) -- (\ref{eq:sparse2}). These fast and slow components
can then be used to form the fast and slow right-hand side functions provided to the
ARKODE multirate integrator.

\begin{figure}[!htb]
    \centering
    \includegraphics[width=0.8\linewidth]{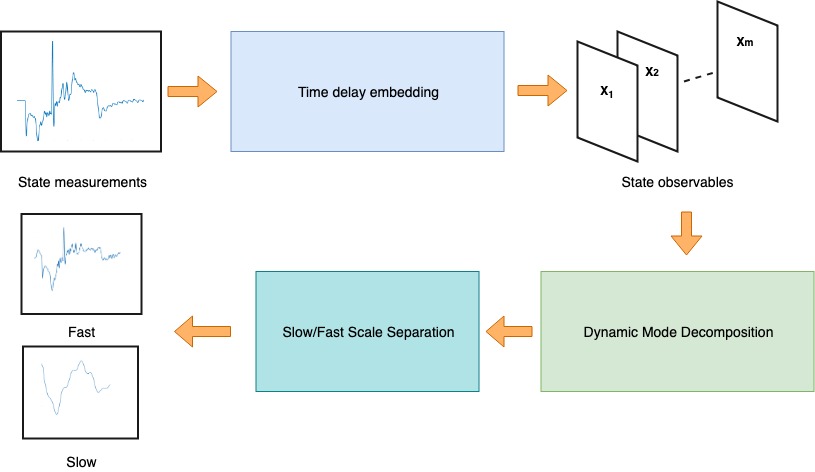}
    \caption{The time-scale separation algrotihm}
    \label{fig:algo}
\end{figure}

\section{Computational Results}\label{s:results}

We applied the time-scale separation algorithm in two different problems. The
first problem is a simple, contrived toy problem. The second scenario uses
synthetic seismic data for tsunami modeling that is closer to real-world data.

\subsection{Toy Problem}

The toy problem uses is defined as the simple IVP:

\[
    \frac{dy}{dt} = \cos(10t) + \sin(t),\quad y(0) = 1.
\]

\begin{figure}[!htb]
    \centering
    \includegraphics[width=0.6\linewidth]{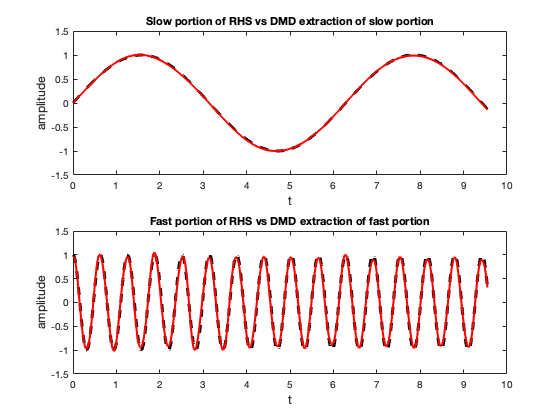}
    \caption{The DMD recreations of the fast and slow components of the
    toy problem right-hand side. These recreations are fed into ARKODE as the
    fast and slow functions.}
    \label{fig:toy_rhs}
\end{figure}

\noindent It is clear the right hand side has a fast component, $\cos(10t)$ and a slow component $\sin(t)$. Using MATLAB we solve the ODE with a uniform time
step size of 0.01. The solution is used as our dynamics data that we 
construct the Hankel matrix with and then perform the DMD. The derivative of the fast
and slow components extracted is shown in Figure~\ref{fig:toy_rhs}. These components
were used in the ARKODE solver by reading the different components in from a data file
instead of computing the values of the functions at some time instance.

\begin{figure}[!htb]
    \centering
    \begin{subfigure}{0.4\linewidth}
        \includegraphics[width=\linewidth]{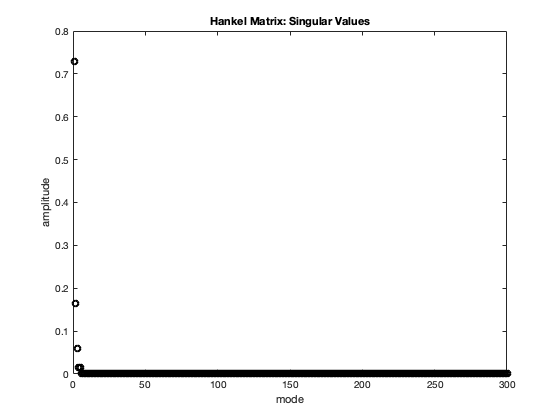}
    \end{subfigure}
    \begin{subfigure}{0.4\linewidth}
        \includegraphics[width=\linewidth]{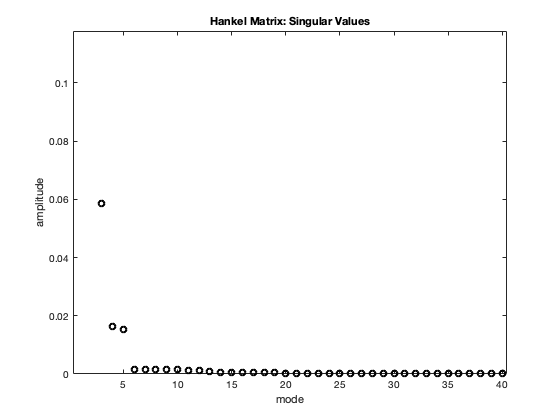}
    \end{subfigure}
    \caption{The singular values of the toy problem $\mathbf{H}$ matrix with 
    300 embeddings. The right image is a zoomed in version of the left. The left 
    image shows the point at which we can truncate when performing the DMD ($\approx20$).}
    \label{fig:toy_hankel}
\end{figure}

\begin{figure}[!htb]
    \centering
    \begin{subfigure}{0.4\linewidth}
        \includegraphics[width=\linewidth]{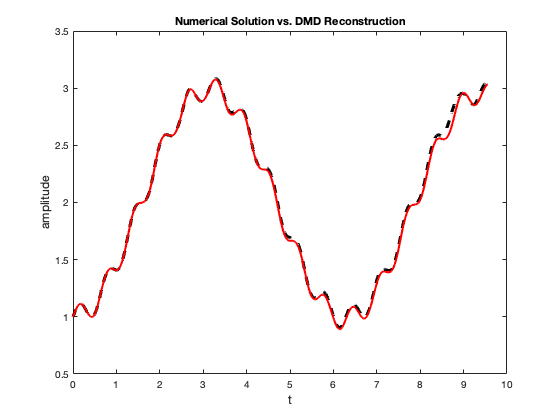}
    \end{subfigure}
    \begin{subfigure}{0.4\linewidth}
        \includegraphics[width=\linewidth]{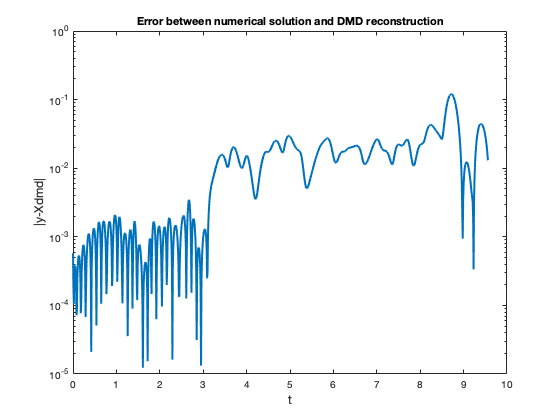}
    \end{subfigure}
    \caption{The image on the left is the DMD reconstruction of the toy
    problem dynamics. The black dotted line is the numerical solution and
    the red line is the DMD approximation. The right image shows the error,
    given by $|\mathbf{X} - \mathbf{X}_\text{DMD}|$, in the approximation.}
    \label{fig:toy_xdmd}
\end{figure}

\subsection{Synthetic Seismic Data}

The second problem we apply our scale separation algorithm to is a seismic problem
where waves, or primary waves, are fast and S waves, or secondary waves are slow.
In \cite{vogl2017high}, Vogl \textit{et al.} present a method for a
high-resolution seismic model of a fault slippage under the ocean. In the
work they present a linear hyperbolic system of equations that model the 2D
plane-strain case of isotropic linear elasticity:

\begin{equation}\label{eq:planestrain}
    q_t + A(x,y)q_x + B(x,y)q_y = 0
\end{equation}

A general Riemann solution to produces the eigenvectors of $n_xA + n_yB$, where
$n = [0~~1]^T$ for the specific data we will examine. The eigenvectors
correspond to P and S waves traveling from the source. 




As part of their work, Vogl \textit{et al.} developed a simulation code for the
problem based on the Clawpack solver package. The simulation code takes
measurements at varying location from the fault with "gauges". The data recorded
corresponds to the vector $q$ which contains key parameters of the plane-strain equations such as the stress, density and velocity. We build a data matrix
$\mathbf{X \in \mathbb{R}^{n \times 10}}$ from 10 sequential gauges and the
vertical velocity. This data matrix is then fed into our time scale separation
algorithm.

\begin{figure}[!htb]
    \centering
    \begin{subfigure}{0.4\linewidth}
        \includegraphics[width=\linewidth]{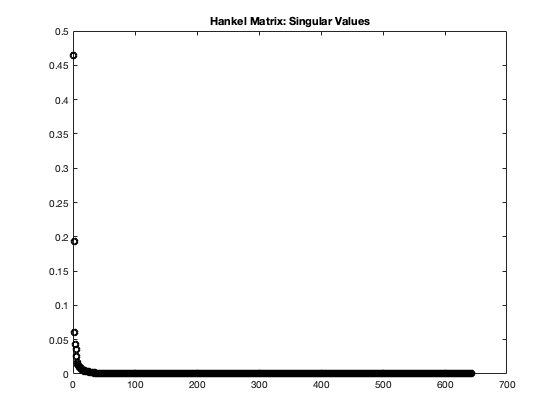}
    \end{subfigure}
    \begin{subfigure}{0.4\linewidth}
        \includegraphics[width=\linewidth]{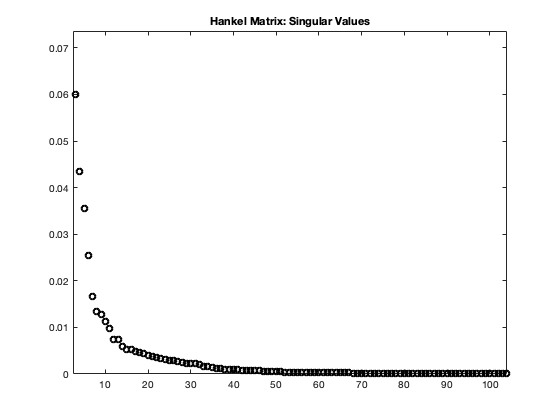}
    \end{subfigure}
    \caption{The singular values of the $\mathbf{H}$ matrix for the seismic data
    with 150 embeddings. The right image is a zoomed in version of the left. The left 
    image shows the point at which we can truncate when performing the DMD ($\approx 75$).}
    \label{fig:toy_hankel}
\end{figure}

\begin{figure}[!htb]
    \centering
    \begin{subfigure}{0.4\linewidth}
        \includegraphics[width=\linewidth]{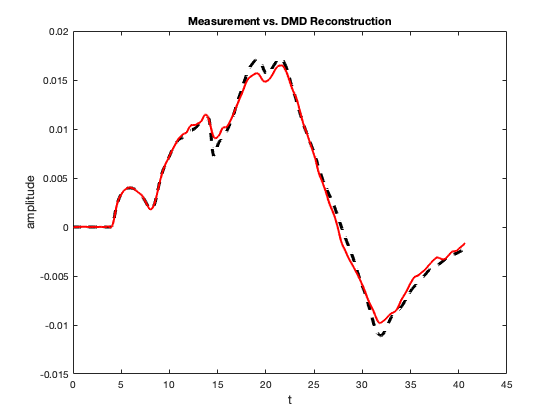}
    \end{subfigure}
    \begin{subfigure}{0.4\linewidth}
        \includegraphics[width=\linewidth]{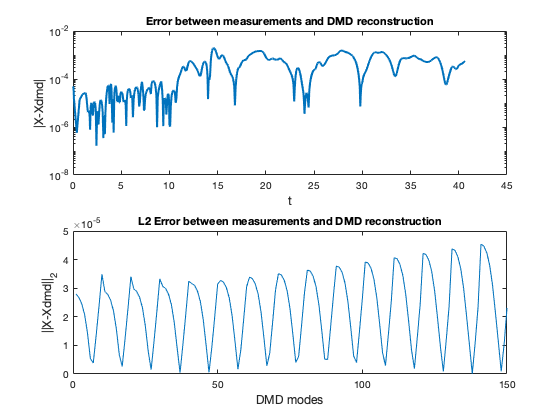}
    \end{subfigure}
    \caption{The image on the left is the DMD reconstruction of the seismic
    problem dynamics. The black dotted line is the numerical solution and
    the red line is the DMD approximation. The right image shows the error,
    given by $|\mathbf{X} - \mathbf{X}_\text{DMD}|$, in the approximation.}
    \label{fig:toy_xdmd}
\end{figure}

\begin{figure}[!htb]
    \centering
    \includegraphics[width=0.6\linewidth]{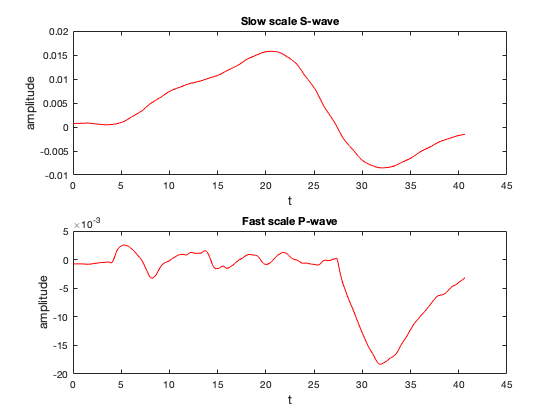}
    \caption{The extracted fast and slow components of the
    vertical velocity in the seismic data.}
    \label{fig:seis_rhs}
\end{figure}

\section{Summary and Conclusions}\label{s:summary}

In summary, in this report, we used dynamic mode decomposition and time delay embedding to extract
the fast and and slow components of the right-hand sides of a simple ODE from data. We then used the extracted components to solve the ODE with ARKODE. Finally, to move towards a real-world use case, we attempted to extract fast and slow scale dynamics from synthetic seismic modeling data. We found
our algorithm to be adequate with the simple dynamics in a toy problem. In the more difficult case
with the synthetic seismic data it showed promise, there seemed to be a split of fast and slow scale
components, but more work needs to be done to verify the accuracy of the components. In the future, we would like to explore using the multi-resolution DMD presented in \cite{kutz2016multiresolution} to compare the performance and fit for the use case.

\section{Acknowledgements}

We would like to thank Christopher Vogl for valuable insight into the seismic model used for the
synthetic data experiment and with assistance in collecting the data from the simulation code.

This work was performed	under the auspices of the U.S.	Department of Energy by	Lawrence Livermore National Laboratory under Contract DE-AC52-07NA27344. LLNL-TR-848209.

\clearpage
\bibliographystyle{siam}
\bibliography{refs.bib}

\section*{Appendix A:\ MATLAB Code}
\begin{lstlisting}[basicstyle=\small,breaklines=true,language=Matlab]
function [Phi,omega,lambda,b,Xdyn,Xdmd] = DMD(X1,X2,r,dt)
% function [Phi,omega,lambda,b,Xdmd] = DMD(X1,X2,r,dt) % Computes the Dynamic Mode Decomposition of X1, X2
%
% INPUTS:
% X1 = X, data matrix
% X2 = X?, shifted data matrix
% Columns of X1 and X2 are state snapshots
% r = target rank of SVD
% dt = time step advancing X1 to X2 (X to X?)
%
% OUTPUTS:
% Phi, the DMD modes
% omega, the continuous-time DMD eigenvalues
% lambda, the discrete-time DMD eigenvalues
% b, a vector of magnitudes of modes Phi
% Xdmd, the data matrix reconstrcted by Phi, omega, b

%% DMD
[U, S, V] = svd(X1, 'econ'); r = min(r, size(U,2));
U_r = U(:, 1:r); % truncate to rank-r
S_r = S(1:r, 1:r);
V_r = V(:, 1:r);
Atilde = U_r' * X2 * V_r / S_r; % low-rank dynamics
[W_r, D] = eig(Atilde);
Phi = X2 * V_r / S_r * W_r;     % DMD modes
lambda = diag(D);               % discrete-time eigenvalues
omega = log(lambda)/dt;         % continuous-time eigenvalues

%% Compute DMD mode amplitudes b
x1 = X1(:, 1);
b = Phi\x1;

%% DMD reconstruction
mm1 = size(X1, 2);
time_dynamics = zeros(r, mm1);
t = (0:mm1-1)*dt; % time vector for iter = 1:mm1,

for iter = 1:mm1
    time_dynamics(:,iter) = (b.*exp(omega*t(iter)));
end

Xdyn = time_dynamics;

Xdmd = Phi * time_dynamics;
\end{lstlisting}
\begin{lstlisting}[basicstyle=\small,breaklines=true,language=Matlab]
close all; clear all; clc

fslow=@(t,y) sin(t);
ffast=@(t,y) cos(10*t);

y0=1;
frhs=@(t,y) [ffast(t,y) + fslow(t,y)];

dt=0.01;
tspan=[0:dt:4*pi];                      % time span
[t,y]=ode45(@(t,y)frhs(t,y),tspan,y0);  % integrate
[~,yslow]=ode45(@(t,y)fslow(t,y),tspan,y0);  % integrate
[~,yfast]=ode45(@(t,y)ffast(t,y),tspan,y0);  % integrate

% time delay embedding
N=300;  % number of embeddings
H=[];
for j=1:N
  H=[H; y(j:length(y)-N+j,:).']; 
end

[U,S,V]=svd(H,'econ'); s=diag(S);

figure(1)
% subplot(2,2,[1 3])
plot(s/sum(s),'ko','Linewidth',2)
title('Hankel Matrix: Singular Values')
xlabel('mode'), ylabel('amplitude')

% subplot(2,2,2)
% plot(U(:,1:5))
% title('Modes')
% subplot(2,2,4)
% plot(V(:,1:5))
% title('Time')


% DMD reconstruction
X1=H(:,1:end-1); X2=H(:,2:end); tt=t(1:size(X1,2));
[Phi,omega,lambda,b,Xpart,Xdmd]=DMD(X1,X2,25,dt);

error=abs(y(1:length(tt))-Xdmd(1,:)');

% plot
figure(2), grid on
plot(tt,X1(1,:),'k--',tt,real(Xdmd(1,:)),'r-','Linewidth',2)
title('Numerical Solution vs. DMD Reconstruction')
xlabel('t','Fontsize',12), ylabel('amplitude','Fontsize',12)

figure(3)
semilogy(tt,error,'Linewidth',2)
title('Error between numerical solution and DMD reconstruction')
xlabel('t','Fontsize',12), ylabel('|y-Xdmd|','Fontsize',12)

%% try and extract fast and slow time-scales

w=abs(omega);
eps=0.15;

Xslow=Phi(:,(w/max(w) < eps))*Xpart((w/max(w) < eps),:);
%Xfast=Phi(:,(w/max(w) >= eps))*Xpart((w/max(w) >= eps),:);
Xfast=Xdmd-abs(Xslow);

figure(4)
subplot(2,1,1)
plot(tt,yslow(1:length(tt)),'k-',tt,Xslow(1,:),'r-')
title('Slow scale portion of the solution')
xlabel('t','Fontsize',12), ylabel('amplitude','Fontsize',12)

subplot(2,1,2)
plot(tt,yfast(1:length(tt)),'k-',tt,Xfast(1,:)+y0,'r-')
title('Fast scale portion of the solution')
xlabel('t','Fontsize',12), ylabel('amplitude','Fontsize',12)


%% get RHS functions back through differentiation

td=tt(1:end-1);

fxslow=diff(Xslow(1,:)')./diff(tt);
fxfast=diff(Xfast(1,:)')./diff(tt);

figure(6)
subplot(2,1,1)
plot(td,fslow(td,y0),'k--',td,fxslow,'r-','Linewidth',2)
title('Slow portion of RHS vs DMD extraction of slow portion')
xlabel('t','Fontsize',12), ylabel('amplitude','Fontsize',12)
subplot(2,1,2)
plot(td,ffast(td,y0),'k--',td,fxfast,'r-','Linewidth',2)
title('Fast portion of RHS vs DMD extraction of fast portion')
xlabel('t','Fontsize',12), ylabel('amplitude','Fontsize',12)


\end{lstlisting}
\begin{lstlisting}[basicstyle=\small,breaklines=true,language=Matlab]
clear all; close all; clc

data1 = dlmread('gauge00030.txt','',3,0);
data2 = dlmread('gauge00031.txt','',3,0);
data3 = dlmread('gauge00032.txt','',3,0);
data4 = dlmread('gauge00033.txt','',3,0);
data5 = dlmread('gauge00034.txt','',3,0);
data6 = dlmread('gauge00035.txt','',3,0);
data7 = dlmread('gauge00036.txt','',3,0);
data8 = dlmread('gauge00037.txt','',3,0);
data9 = dlmread('gauge00038.txt','',3,0);
data10 = dlmread('gauge00039.txt','',3,0);

% extract vertical? velocity
j = 7;
t = data1(:,2);
p = [data1(:,j), data2(:,j), data3(:,j), data4(:,j),...
     data5(:,j), data6(:,j), data7(:,j), data8(:,j),...
     data9(:,j), data10(:,j)];

figure(1)
plot(t,p(:,1),'Linewidth',2)

% V^-1 qt + WV^-1 qx

%% time delay embedding

N=150;  % number of embeddings
H=[];
for j=1:N
  H=[H; p(j:length(p)-N+j,:).']; 
end


[u,sigma,v] = svd(H,'econ'); s = diag(sigma);

figure(2)
plot(s/sum(s),'ko','Linewidth',2)
title('Hankel Matrix: Singular Values')

% figure(3)
% subplot(2,1,1)
% plot(u(:,1:N))
% title('Modes')
% subplot(2,1,2)
% plot(v(:,1:N))
% title('Time')

%% DMD reconstruction

X1=H(:,1:end-1); X2=H(:,2:end); dt=t(2)-t(1); tt=t(1:size(X1,2));
[Phi,omega,lambda,b,Xdyn,Xdmd]=DMD(X1,X2,75,dt);

error=abs(X1(1,:)-Xdmd(1,:));
l2error=(1/size(X1,2))*sqrt(abs(sum((X1(1:N,:)'-Xdmd(1:N,:)')).^2));

% plot
figure(4), grid on
plot(tt,X1(1,:),'k--',tt,real(Xdmd(1,:)),'r-','Linewidth',2)
title('Measurement vs. DMD Reconstruction')
xlabel('t','Fontsize',12), ylabel('amplitude','Fontsize',12)

figure(5)
plot(tt,real(Xdyn),'Linewidth',2)
title('DMD Dynamics'), xlabel('time'), ylabel('amplitude')

figure(6)
plot(real(omega),imag(omega),'ko','Linewidth',2)
title('Eigenvalues'), xlabel('real'), ylabel('imaginary')

figure(7)
subplot(2,1,1)
semilogy(tt,error,'Linewidth',2)
title('Error between measurements and DMD reconstruction')
xlabel('t','Fontsize',12), ylabel('|X-Xdmd|','Fontsize',12)
subplot(2,1,2)
plot(1:N,l2error)
title('L2 Error between measurements and DMD reconstruction')
xlabel('DMD modes','Fontsize',12), ylabel('||X-Xdmd||_2','Fontsize',12)

%% try and extract fast and slow time-scales

w=abs(omega);
eps=0.05;

Xslow=Phi(:,(w/max(w) < eps))*Xdyn((w/max(w) < eps),:);
%Xfast=Phi(:,(w/max(w) >= eps))*Xpart((w/max(w) >= eps),:);
Xfast=Xdmd-abs(Xslow);

figure(8)
subplot(2,1,1)
plot(tt,real(Xslow(1,:)),'r-')
title('Slow scale S-wave')
xlabel('t','Fontsize',12), ylabel('amplitude','Fontsize',12)
subplot(2,1,2)
plot(tt,real(Xfast(1,:)),'r-')
title('Fast scale P-wave')
xlabel('t','Fontsize',12), ylabel('amplitude','Fontsize',12)

figure(9), plot(tt,real(Xslow(1,:)),'Linewidth',2)
figure(10), plot(tt,real(Xfast(1,:)),'Linewidth',2)


\end{lstlisting}
 
\end{document}